\documentclass[11pt,leqno]{article}

\usepackage{amsmath,amsfonts,amscd,amssymb,theorem}

\long\def\comment#1\endcomment{}

\makeatletter
\begingroup
\gdef\th@dotted{\normalfont\itshape
  \def\@begintheorem##1##2{%
        \item[\hskip\labelsep \theorem@headerfont ##1\ ##2.]}%
\def\@opargbegintheorem##1##2##3{%
   \item[\hskip\labelsep \theorem@headerfont ##1\ ##2\ (##3).]}}
\endgroup
\makeatother

\theoremstyle{dotted}

\newtheorem{theorem}{Theorem}[section]
\newtheorem{lemma}[theorem]{Lemma}

\newtheorem{prop}[theorem]{Proposition}
\newtheorem{corr}[theorem]{Corollary}

\makeatletter
\begingroup
\gdef\th@upshape{\normalfont
  \def\@begintheorem##1##2{%
        \item[\hskip\labelsep \theorem@headerfont ##1\ ##2.]}%
\def\@opargbegintheorem##1##2##3{%
   \item[\hskip\labelsep \theorem@headerfont ##1\ ##2\ (##3).]}}
\endgroup
\makeatother

\theoremstyle{upshape}

\newtheorem{defn}[theorem]{Definition}
\newtheorem{remark}[theorem]{Remark}
\newtheorem{example}[theorem]{Example}

\makeatletter
\renewcommand{\subsection}{\@startsection{subsection}{2}{0pt}{-3ex
plus -1ex minus -0.2ex}{-2mm plus -0pt minus
-2pt}{\normalfont\bfseries}} \makeatother

\makeatletter
\@addtoreset{equation}{section}
\makeatother

\newcommand{\cntrct}                
{\hspace{2pt}\raisebox{1pt}{\text{$\lrcorner$}}\hspace{2pt}}

\newcommand{\proof}[1][Proof.]{\smallskip\noindent{\em #1}}
\def\endproof{\hfill\ensuremath{\square}\par\medskip}

\def\eqref#1{\thetag{\ref{#1}}}

\let\latexref=\ref
\def\ref#1{{\normalfont{\latexref{#1}}}}

\newcommand{\wt}{\widetilde}
\newcommand{\wh}{\widehat}

\newcommand{\idot}{{\:\raisebox{1pt}{\text{\circle*{1.5}}}}}
\newcommand{\hdot}{{\:\raisebox{3pt}{\text{\circle*{1.5}}}}}
\newcommand{\Spec}{\operatorname{Spec}}

\def\dlim_#1{{\displaystyle\lim_{#1}}^\hdot}

\newcommand{\Vect}{\operatorname{\it\!-Vect}}
\newcommand{\Sets}{\operatorname{Sets}}
\newcommand{\bimod}{\operatorname{\!-\sf bimod}}

\newcommand{\cchar}{\operatorname{\sf char}}

\newcommand{\Fun}{\operatorname{Fun}}

\newcommand{\N}{{\mathbb N}}

\newcommand{\Sec}{\operatorname{\sf Sec}}

\newcommand{\defeq}{\overset{\text{\sf\tiny def}}{=}}

\newcommand{\gr}{\operatorname{\sf gr}}

\newcommand{\C}{{\mathcal C}}

\newcommand{\B}{{\mathcal B}}

\newcommand{\D}{{\mathcal D}}

\newcommand{\Ext}{\operatorname{Ext}}

\newcommand{\I}{\operatorname{{\sf I}}}

\newcommand{\End}{\operatorname{End}}

\newcommand{\ppt}{\operatorname{{\sf pt}}} 

\newcommand{\Ker}{\operatorname{{\sf Ker}}}
\newcommand{\Coker}{\operatorname{{\sf Coker}}}

\newcommand{\id}{\operatorname{\sf id}}

\newcommand{\tr}{\operatorname{\sf tr}}

\newcommand{\Tor}{\operatorname{Tor}}
\newcommand{\Shv}{\operatorname{Shv}}

\newcommand{\m}{{\mathfrak m}}

\title{Cyclic homology with coefficients}

\author{D. Kaledin\thanks{Partially supported by CRDF grant
RUM1-2694-MO05.}}

\date{\em To Yu.~I.~Manin, the founder, on the occasion of his 70-th
birthday}

\begin{document}

\maketitle

\tableofcontents

\section*{Introduction}

Ever since it was discovered in 1982 by A. Connes \cite{C1} and
B. Tsygan \cite{tsy}, cyclic homology occupies a strange place in
the realm of homological algebra. Normally in homological algebra
problems, one expects to start from some data, such e.g. a
topological space $X$, then construct some abelian category, such as
the category of sheaves on $X$, and then define the cohomology of
$X$ by computing the derived functors of some natural functor, such
as e.g. the global sections functor $\Gamma(X,-)$. Admittedly, this
is a modern formulation, but it had certainly been current already
in 1982. Cyclic homology starts with an associative algebra $A$, and
defines its homology groups $HC_\idot(A)$, but there are absolutely
no derived functors in sight. Originally, $HC_\idot(A)$ were defined
as the homology of an explicit complex, -- which anyone trained to
use triangulated categories cannot help but take as an insult. Later
A. Connes \cite{C} improved on the definition by introducing the
abelian category of so-called cyclic vector spaces. However, the
passage from $A$ to its associated cyclic vector space $A_\#$ is
still done by an explicit {\em ad hoc} formula. It is as if we were
to know the bar-complex which computes the homology of a group,
without knowing the definition of the homology of a group.

This situation undoubtedly irked many people over the years, but to
the best of my knowledge, no satisfactory solution has been
proposed, and it may not exist -- indeed, many relations to the de
Rham homology notwithstanding, it is not clear whether cyclic
homology properly forms a part of homological algebra at all (to the
point that e.g. in \cite{FT} the word ``homology'' is not used at
all for $HC_\idot(A)$, and it is called instead {\em additive
$K$-theory} of $A$). In the great codification of homological
algebra done in \cite{GM1}, cyclic homology only appears in the
exercises. This is not surprising, since the main unifying idea of
\cite{GM1} is the ideology of ``linearization'': homological algebra
linearizes geometry, just as functional analysis used to do 50 years
ago; triangulated categories and adjoint functors are modern-day
versions of Banach spaces and adjoint linear operators. This has
been an immensely successful and clarifying point of view, in
general, but $HC_\idot(A)$ sticks out on a complete tangent -- there
is simply no natural place for it in this framework.

\medskip

This paper arouse as one more attempt to propose a solution to the
difficulty -- to find a natural triangulated category where
$HC_\idot(-)$ would be able to live with a certain level of comfort
(and with all the standard corollaries such as the notion of cyclic
homology with coefficients, the ability to compute cyclic homology
by whatever resolution is convenient, not just the bar resolution,
and so on).

\medskip

In a sense, our attempt has been successful: we define a
triangulated category which can serve as the natural ``category of
coefficients'' for cyclic homology of an algebra $A$, and we prove
the comparison theorem that shows that when the coefficients are
trivial, the new definition of cyclic homology is equivalent to the
old one. In fact, the algebra $A$ enters into the construction only
through the category $A\bimod$ of $A$-bimodules; we also show how to
generalize the construction so that $A\bimod$ is replaced with a
more general tensor abelian category $\C$.

From a different point of view, though, out attempt failed
miserably: the correspondence $A \mapsto A_\#$, being thrown out of
the window, immediately returns through the door in a new and
``higher-level disguise'': it is now applied not to the algebra $A$,
but to the tensor category $\C = A\bimod$. Then in practice, the
freedom to choose an arbitrary resolution to compute the derived
functors leads, in our approach to $HC_\idot(-)$, to complexes which
are even larger than the original complex, and at some point the
whole exercise starts to look pointless.

Still, we believe that, all said and done, some point can be found,
and some things are clarified in our approach; one such thing is,
for instance, the version of Gauss-Manin connection for cyclic
homology discovered by E. Getzler \cite{getz}. Beside, we do propose
a definition of cyclic homology which makes sense for a general
tensor category; and in some particular questions, even the
computations can be simplified. As for the presence of the
$A_\#$-construction, this might be in the nature of things, after
all -- not a bug of the theory, but a necessary feature. However, we
best leave it to the reader to be the judge.

\medskip

The paper is organized as follows. In Section 1 we recall A. Connes'
second definition of cyclic homology which uses the cyclic category
$\Lambda$; we also recall some facts about homology of small
categories that we will need. We have tried to give only the
absolute minimum -- the reader not familiar with the material will
have to consult the references. In Section 2 we introduce our main
object: the notion of a {\em cyclic bimodule} over an associative
algebra $A$, and the derived category of such bimodules. We also
introduce cyclic homology $HC_\idot(A,M)$ with coefficients in a
cyclic bimodule $M$. In Section 3 we give a very short derivation of
the Gauss-Manin connection; strictly speaking, the language of
cyclic bimodules is not needed for this, but we believe that it
shows more clearly what is really going on. In Section 4, we show
how to replace the category $A\bimod$ everywhere with a more general
tensor abelian category $\C$. Section 5 is a postface, or a
``discussion'' (as they do in medical journals) -- we discuss some
of the further things one might (and should) do with cyclic
bimodules, and how to correct some deficiencies of the theory
developed in Sections 2 and 4.

\subsection*{Acknowledgements.} In the course of this work, I have
benefited greatly from discussions with A. Beilinson, E. Getzler,
V. Ginzburg, A. Kuznetsov, N. Markarian, D. Tamarkin, and
B. Tsygan. I am grateful to Northwestern Univeristy, where part of
this work was done, and where some of the results were presented in
seminars, with great indulgence from the audience towards the
unfinished state they were in. And, last but not least, it is a
great pleasure and a great opportunity to dedicate the paper to Yuri
Ivanovich Manin on his birthday. Besides all the usual things, I
would like to stress that it is the book \cite{GM1}, -- and
\cite{GM2}, to a lesser extent -- which shaped the way we look at
homological algebra today, at least ``we'' of my generation and of
Moscow school. Without Manin's decisive influence, this paper
certainly would not have appeared (as in fact at least a half of the
papers I ever wrote).

\section{Recollection on cyclic homology.}

We start by recalling, extremely briefly, A. Connes' approach to
cyclic homology, which was originally introduced in \cite{C} (for
detailed overviews, see e.g. \cite[Section 6]{L} or
\cite[Appendix]{FT}; a brief but complete exposition
using the same language and notation as in this paper can be found
in \cite[Section 1]{Ka}).

Connes' approach relies on the technique of homology of small
categories. Fix a base field $k$. Recall that for every small
category $\Gamma$, the category $\Fun(\Gamma,k)$ of functors from
$\Gamma$ to the category $k\Vect$ of $k$-vector spaces is an abelian
category with enough projectives and enough injectives, with derived
category $\D(\Gamma,k)$. For any object $E \in \Fun(\Gamma,k)$, the
homology $H_\idot(\Gamma,E)$ of the category $\Gamma$ with
coefficients in $E$ is by definition the derived functor of the
direct limit functor
$$
\displaystyle\lim_{\overset{\to}{\Gamma}}:\Fun(\Gamma,k) \to k\Vect.
$$
Analogously, the cohomology $H^\hdot(\Lambda,E)$ is the derived
functor of the inverse limit
$\displaystyle\lim_{\overset{\gets}{\Gamma}}$. Equivalently,
$$
H^\hdot(\Gamma,E) = \Ext^\hdot(k,E),
$$
where $k \in \Fun(\Gamma,k)$ is the constant functor (all objects in
$\Gamma$ go to $k$, all maps go to identity). In particular,
$H^\hdot(\Gamma,k)$ is an algebra. For any $E \in \Fun(\Gamma,k)$,
the cohomology $H^\hdot(\Gamma,E)$ and the homology
$H_\idot(\Gamma,E)$ are modules over $H^\hdot(\Gamma,k)$.

We also note, although it is not needed for the definition of cyclic
homology, that for any functor $\gamma:\Gamma' \to \Gamma$ between
two small categories, we have the pullback functor
$\gamma^*:\Fun(\Gamma,k) \to \Fun(\Gamma',k)$, and for any $E \in
\Fun(\Gamma,k)$, we have natural maps
\begin{equation}\label{dir.im}
H_\idot(\Gamma',\gamma^*E) \to H_\idot(\Gamma,E),\qquad
H^\hdot(\Gamma,E) \to H^\hdot(\Gamma',\gamma^*E).
\end{equation}
Moreover, the pullback functor $\gamma^*$ has a left adjoint
$\gamma_!:\Fun(\Gamma',k) \to \Fun(\Gamma,k)$ and a right-adjoint
$f_*:\Fun(\Gamma',k) \to \Fun(\Gamma,k)$, known as the left and
right Kan extensions. In general, $f_!$ is right-exact but it need
not be left-exact. We will need one particular case where it is
exact. Assume given a covariant functor $V:\Gamma \to \Sets$ from a
small category $\Gamma$ to the category of sets, and consider the
category $\Gamma'$ of pairs $\langle [a],v \rangle$ of an object
$[a] \in \Gamma$ and an element $v \in V([a])$ (maps in $\Gamma'$
are those maps $\gamma:[a] \to [a']$ which send $v \in V([a])$ to
$v' \in V([a'])$. Such a category is known as a {\em discrete
cofibration} over $\Gamma$ associated to $V$, see \cite{SGA}. Then
the Kan extension $f_!$ associated to the forgetful functor
$f:\Gamma' \to \Gamma$, $\langle [a],v \rangle \mapsto [a]$ is
exact, and is easy to compute: for any $E \in \Fun(\Gamma',k)$ and
$[a] \in \Gamma$, we have
\begin{equation}\label{discr}
f_!E([a]) = \bigoplus_{v \in V([a])}E(\langle [a],v \rangle).
\end{equation}
Moreover, for any $E \in \Fun(\Gamma,k)$, this imediately gives the
projection formula:
\begin{equation}\label{projj}
f_!f^*E \cong E \otimes F_!k,
\end{equation}
where, as before, $k \in \Fun(\Gamma',k)$ stands for the constant
functor.

\medskip

For applications to cyclic homology, one starts with introducing the
{\em cyclic category} $\Lambda$. This is a small category whose
objects $[n]$ are numbered by positive integers $n \geq 1$. One
thinks of an object $[n]$ as a circle $S^1$ with $n$ distinct marked
points; we denote the set of these points by $V([n])$. The set of
maps $\Lambda([n'],[n])$ from $[n']$ to $[n]$ is then the set of
homotopy classes of continuous maps $f:S^1 \to S^1$ such that
\begin{itemize}
\item $f$ has degree $1$, sends marked points to marked points, and
  is non-decreasing with respect to the natural cyclic order on
  $S^1$ (that is, if a point $a \in S^1$ lies between points $b$ and
  $c$ when counting clockwise, then the same is true for $f(a)$,
  $f(b)$ and $f(c)$).
\end{itemize}
In particular, we have $\Lambda([1],[n]) = V([n])$. This topological
description of the cyclic category $\Lambda$ is easy to visualize,
but there are also alternative combinatorial descriptions
(e.g. \cite[Exercize II.1.6]{GM1}, \cite[Section 6]{L}, or
\cite[A.2]{FT}, retold in \cite[Section 1.4]{Ka}). All the
descriptions are equivalent. Objects in $\Fun(\Lambda,k)$ are
usually called {\em cyclic vector spaces}.

The cyclic category $\Lambda$ is related to the more familiar {\em
simplicial category} $\Delta^{opp}$, the opposite to the category
$\Delta$ of finite non-empty linearly ordered sets. To understand
the relation, consider the discrete cofibration
$\Lambda_{[1]}/\Lambda$ associated to the functor $V:\Lambda \to
\Sets$ -- equivalently, $\Lambda_{[1]}$ is the category of objects
$[n]$ in $\Lambda$ eqipped with a map $[1] \to [n]$. Then it is easy
to check that $\Lambda_{[1]}$ is equivalent to the
$\Delta^{opp}$. From now on, we will abuse the notation and identify
$\Lambda_{[1]}$ and $\Delta^{opp}$. We then have a natural
projection $\Delta^{opp} = \Lambda_{[1]} \to \Lambda$, $\langle
[n],v \rangle \mapsto [n]$, which we denote by $j:\Delta^{opp} \to
\Lambda$.

For any cyclic $k$-vector space $E \in \Fun(\Lambda,k)$, we have its
restriction $j^*E \in \Fun(\Delta^{opp},E)$, a simplicial vector
space. One defines the cyclic homology $HC_\idot(E)$ and the
Hochschild homology $HH_\idot$ of $E$ by
$$
HC_\idot(E) \defeq H_\idot(\Lambda,E), \qquad HH_\idot(E) \defeq
H_\idot(\Delta^{opp},j^*E).
$$
By \eqref{dir.im}, we have a natural map $HH_\idot(E) \to
HC_\idot(E)$ (moreover, since $j:\Delta^{opp} \to \Lambda$ is a
discrete cofibration, the Kan extension $j_!$ is exact, so that we
have $HH_\idot(E) \cong HC_\idot(j_!j^*E)$, and the natural map is
induced by the adjunction map $j_!j^*E \to E$).  It has been shown
by A. Connes that this map fits into a long exact sequence
\begin{equation}\label{connes}
\begin{CD}
HH_\idot(E) @>>> HC_\idot(E) @>{u}>> HC_{\idot-2}(E) @>>>.
\end{CD}
\end{equation}
Here the map $u$ is the so-called {\em periodicity map} on
$HC_\idot(E)$: one shows that the algebra $H^\hdot(\Lambda,k)$ is
isomorphic to the polynomial algebra $k[u]$ in one generator $u$ of
degree $2$, and the periodicity map on homology is simply the action
of this generator. This allows to define a third homological
invariant, the {\em periodic cyclic homology} $HP_\idot(E)$ -- to do
it, one inverts the periodicity map.

\begin{defn}
For any cyclic $k$-vector space $E \in \Fun(\Lambda,k)$, the {\em
periodic cyclic homology} of $E$ is defined by
$$
HP_\idot(E) = \dlim_{\overset{u}{\gets}}HC_\idot(E),
$$
where $\dlim_{\gets}$ denotes the derived functor of the
inverse limit $\displaystyle\lim_{\gets}$.
\end{defn}

Assume now given an associative unital algebra $A$ over $k$. To
define its cyclic homology, we associate to $A$ a canonical cyclic
vector space $A_\#$ in the following way. We set $A_{\#}([n]) =
A^{\otimes V([n])}$, the tensor product of $n$ copies of the vector
space $A$ numbered by marked points $v \in V([n])$. Then for any map
$f \in \Lambda([n'],[n])$, we define
\begin{equation}\label{hash.def}
A_\#(f) = \bigotimes_{v \in V([n])}m_{f^{-1}(v)}:A^{\otimes V([n'])} =
\bigotimes_{v \in V([n])}A^{\otimes f^{-1}(v)} \to A^{\otimes
  V([n])},
\end{equation}
where for any linearly ordered finite set $S$, $m_S:A^{\otimes S}
\to A$ is the canonical multiplication map induced by the
associative algebra structure on $A$ (and if $S$ is empty, we set
$A^{\otimes S} = k$, and $m_S$ is the embedding of the unity). This
is obviously compatible with compositions, and it is well-defined
since for any $v \in V([n])$, its preimage $f^{-1} \subset V([m])$
carries a natural linear order induced by the orientation of the
circle $S^1$.

\begin{defn}\label{alg.def}
For any associative unital algebra $A$ over $k$, its Hochschild,
cyclic and periodic cyclic homology $HH_\idot(A)$, $HC_\idot(A)$,
$HP_\idot(A)$ is defined as the corresponding homology of the cyclic
$k$-vector space $A_\#$:
$$
HH_\idot(A) \defeq HH_\idot(A_\#),\quad
HC_\idot(A) \defeq HC_\idot(A_\#),\quad
HC_\idot(P) \defeq HP_\idot(A_\#).
$$
\end{defn}

\section{Cyclic bimodules.}\label{naive}

Among all the homology functors introduced in
Definition~\ref{alg.def}, Hochschild homology is the most accesible,
and this is because it has another definition: for any associative
unital algebra $A$ over $k$, we have
\begin{equation}\label{hh.def}
HH_{\idot} = \Tor^\hdot_{A^{opp} \otimes A}(A,A),
\end{equation}
where $\Tor^\hdot$ is taken over the algebra $A^{opp} \otimes A$
(here $A^{opp}$ denotes $A$ with the multiplication taken in the
opposite direction).

This has a version with coefficients: if $M$ is a left module over
$A^{opp} \otimes A$, -- in other words, an $A$-bimodule, -- one
defines Hochschild homology of $A$ with coefficients in $M$ by
\begin{equation}\label{hoch.coeff}
HH_\idot(A,M) = \Tor^\hdot_{A^{opp} \otimes A}(M,A).
\end{equation}
The category $A\bimod$ of $A$-bimodules is a unital (non-symmetric)
tensor category, with tensor product $- \otimes_A -$ and the unit
object $A$. Hochschild homology is a homological functor from
$A\bimod$ to $k\Vect$.

To obtain a small category interpretation of $HH_\idot(A,M)$, one
notes that for any $n,n' \geq 0$, the $A$-bimodule structure on $M$
induces a multiplication map
$$
A^{\otimes n} \otimes M \otimes A^{\otimes n'} \to M.
$$
Therefore, if to any $\langle [n],v \rangle \in \Delta^{opp}$ we
associate the $k$-vector space 
\begin{equation}\label{M.Delta}
M^\Delta_\#([n]) = M \otimes A^{\otimes (V([n]) \setminus \{v\})},
\end{equation}
with $M$ filling the place corresponding to $v \in V([n])$, then
\eqref{hash.def} make perfect sense for those maps $f:[n'] \to [n]$
which preserve the distinguished points. Thus to any $M \in
A\bimod$, we can associate a simplicial $k$-vector space
$M^\Delta_\# \in \Fun(\Delta^{opp},k)$. In the particular case
$M=A$, we have $A_\#^\Delta = j^*A_\#$.

\begin{lemma}\label{hoch}
For any $M \in A\bimod$, we have a canonical isomorphism
\begin{equation}\label{hh.iso}
HH_\idot(A,M) \cong H_\idot(\Delta^{opp},M^\Delta_\#).
\end{equation}
\end{lemma}

\proof{} It is well-known that for any simplicial $k$-vector space
$E$, the homology $H_\idot(\Delta^{opp},E)$ can be computed by the
standard complex of $E$ (that is, the complex with terms $E([n])$
and the differential $d = \sum_i(-1)^id_i$, where $d_i$ are the face
maps). In particular, $H_0(\Delta^{opp},M^\Delta_\#)$ is the
cokernel of the map $d:A \otimes M \to M$ given by $d(a \otimes m) =
am-ma$. The natural projection $M \to M \otimes_{A^{opp} \otimes A}
A$ obviously factors through this cokernel, so that we have a
natural map
$$
\rho_0:H_0(\Delta^{opp},M^\Delta_\#) \to HH_0(A,M).
$$
Both sides of \eqref{hh.iso} are homological functors in $M$, and
$HH_\idot(A,M)$ is a universal homological functor (=the derived
functor of $HH_0(A,M)$); therefore the map $\rho_0$ extends to a map
$\rho_\idot:H_\idot(\Delta^{opp},M^\Delta_\#) \to HH_\idot(A,M)$.
To prove that $\rho_\hdot$ is an isomorphism for any $M$, it
suffices to prove it when $M$ is free over $A^{opp} \otimes A$, or
in fact, when $M = A^{opp} \otimes A$. Then on one hand, $HH_0(A,M)
= A$, and $HH_i(A,M) = 0$ for $i \geq 1$. And on the other hand, the
standard complex associated to the simplicial $k$-vector space
$(A^{opp} \otimes A)_\#^\Delta$ is just the usual bar resolution of
the diagonal $A$-bimodule $A$.
\endproof

It is more or less obvious that for an arbitrary $M \in A\bimod$,
$M_\#^\Delta$ does not extend to a cyclic vector space -- in order
to be able to define $HC_\idot(A,M)$, we have to equip the bimodule
$M$ with some additional structure. To do this, we want to use the
tensor structure on $A\bimod$. The slogan is the following:
\begin{itemize}
\item To find a suitable category of coefficients for cyclic
homology, we have to repeat the definition of the cyclic vector
space $A_\# \in \Fun(\Lambda,k)$, but replace the associative
algebra $A$ in this definition with the tensor category $A\bimod$.
\end{itemize}
Let us explain what this means.

First, consider an arbitrary associative unital monoidal category
$\C$ with unit object $I$ (at this point, not necessarily
abelian). For any integer $n$, we have the Cartesian product $\C^n =
\C \times \C \times \dots \times \C$. Moreover, the product on $\C$
induces a product functor
$$
m:\C^n \to \C,
$$
where if $n=0$, we let $\C^n=\ppt$, the category with one object and
one morphism, and let $m:\ppt \to \C$ be the embedding of the unit
object. More generally, for any finite linearly ordered set $S$ with
$n$ elements, we have a product functor $m_S:\C^S \to \C$, where
$\C^S = \C^n$ with multiples in the product labeled by elements of
$S$. Then for any $[n],[n'] \in \Lambda$, and any $f:[n'] \to [n]$,
we can define a functor $f_!:\C^{V([n'])} \to \C^{V([n])}$ by the
same formula as in \eqref{hash.def}:
\begin{equation}\label{trans.f}
f_! = \prod_{v \in V([n])}m_{f^{-1}(v)}:\C^{V([n'])} = \prod_{v \in
  V([n])}\C^{f^{-1}(v)} \to \C^{V([n])}.
\end{equation}
The natural associativity isomorphism for the product on $\C$
induces natural isomorphisms $(f \circ f')_! \cong f_!  \circ f'_!$,
and one checks easily that they satisfy natural compatibility
conditions. All in all, setting $[n] \mapsto \C^{V([n])}$, $f
\mapsto f_!$ defines a weak functor (a.k.a.\ lax functor, a.k.a.\
$2$-fun\-c\-tor, a.k.a.\ pseudofunctor in the original terminology
of Grothendieck) from $\Lambda$ to the category of categories.
Informally, we have a ``cyclic category''.

To work with weak functors, it is convenient to follow
Grothendieck's approach in \cite{SGA}. Namely, instead of
considering a weak functor directly, we define a {\em category}
$\C_\#$ in the following way: its objects are pairs $\langle [n],M_n
\rangle$ of an object $[n]$ of $\Lambda$ and an object $M_n \in
\C^n$, and morphisms from $\langle [n'],M_{n'} \rangle$ to $\langle
[n],M_n\rangle$ are pairs $\langle f,\iota_f\rangle$ of a map
$f:[n'] \to [n]$ and a bimodule map $\iota_f:f_!(M_{n'}) \to M_n$. A
map $\langle f, \iota_f \rangle$ is called {\em cocartesian} if
$\iota_f$ is an isomorphism. For the details of this construction,
-- in particular, for the definition of the composition of
morphisms, -- we refer the reader to \cite{SGA}.

The category $\C_\#$ comes equipped with a natural forgetful
projection $\tau:\C_\# \to \Lambda$, and this projection is a {\em
cofibration} in the sense of \cite{SGA}. A {\em section} of this
projection is a functor $\sigma:\Lambda \to \C_\#$ such that $\tau
\circ \sigma = \id$ (since $\Lambda$ is small, there is no harm in
requiring that two functors from $\Lambda$ to itself are equal, not
just isomorphic). These sections obviously form a category which we
denote by $\Sec(\C_\#)$. Explicitly, an object $M_\# \in
\Sec(\C_\#)$ is given by the following:
\begin{enumerate}
\item a collection of objects $M_n = M_\#([n]) \in \C^n$, and
\item a collection of transition maps $\iota_f:f_!M_{n'} \to M_n$
  for any $n$, $n'$, and $f \in \Lambda([n'],[n])$,
\end{enumerate}
subject to natural compatibility conditions. 

A section $\sigma:\Lambda \to \C_\#$ is called cocartesian if
$\sigma(f)$ is a cocartesian map for any $[n],[n'] \in \Lambda$ and
$f:[n'] \to [n]$ -- equivalently, a section is cocartesian if all
the transition maps $\iota_f$ are isomorphisms. Cocartesian sections
form a full subcategory $\Sec_{cart}(\C_\#)$

\begin{lemma}\label{cycl.str}
The category $\Sec_{cart}(\C_\#)$ of cocartesian objects $M_\# \in
\Sec(\C_\#)$ is equivalent to the category of the following data:
\begin{enumerate}
\item an object $M = M_\#([1]) \in \C$, and
\item an isomorphism $\tau: I \times M \to M \times I$ in the
  category $\C^2 = \C \times \C$,
\end{enumerate}
such that, if we denote by $\tau_{ij}$ the endomorphism of $I \times
I \times M \in \C^3$ obtained by applying $\tau$ to the $i$-th and
$j$-th multiple, we have $\tau_{31} \circ \tau_{12} \circ \tau_{23}
= \id$.
\end{lemma}

\proof{} Straghtforward and left to the reader. \endproof

Thus the natural forgetfull functor $\Sec_{cart}(\C_\#) \to \C$,
$M_\# \mapsto M_\#([1])$ is faithful: an object in
$\Sec_{cart}(\C_\#)$ is given by $M_\#([1])$ plus some extra
structure on it, and all the higher components $M_\#([n])$, $n \geq
2$, together with the transition maps $\iota_f$, can be recovered
from $M_\#([1])$ and this extra structure.

Return now to the abelian situation: we are given an associative
unital algebra $A$ over a field $k$, and our monoidal category is
$\C = A\bimod$, with the natural tensor product. Then for every $n$,
the product $A\bimod^n$ has a fully faithful embedding $A\bimod^n
\to A^{\otimes n}\bimod$, $M_1 \times M_2 \times \dots \times M_n
\mapsto M_1 \boxtimes M_2 \boxtimes \dots \boxtimes M_n$, and one
checks easily that the multiplication functors $m_S$ actually extend
to right-exact functors
$$
m_S:A^{\otimes S}\bimod \to A\bimod;
$$
for instance, one can define $m_S$ as
$$
m_S(M) = M/\{ a_{v'}m - ma_v \mid v \in S, a \in A, m \in M \},
$$
where $a_v = 1 \otimes \dots \otimes a \otimes \dots \otimes 1 \in
A^{\otimes S}$ with $a$ at the $v$-th position, and $v' \in S$ is
the next element after $v$. We can therefore define the cofibered category
$A\bimod_\#/\Lambda$ with fiber $A^{\otimes V([n])}\bimod$ over $[n]
\in \Lambda$, and transition functors $f_!$ as in
\eqref{trans.f}. We also have the category of sections
$\Sec(A\bimod_\#)$ and the subcategory of cocartesian sections
$\Sec_{cart}(A\bimod_\#) \subset \Sec(A\bimod_\#)$.

\begin{lemma}\label{sec.ab}
The category $\Sec(A\bimod_\#)$ is a $k$-linear abelian category.
\end{lemma}

\proof[Sketch of a proof.] This is a general fact about cofibered
categories; the proof is straightforward. The kernel $\Ker\phi$ and
cokernel $\Coker\phi$ of a map $\phi:M_\# \to M'_\#$ between objects
$M_\#,M'_\# \in \Sec(A\bimod_\#)$ are taken pointwise: for every
$n$, we have an exact sequence
$$
0 \to (\Ker\phi)([n]) \to M_\#([n]) \overset{\phi}{\to} 
M'_\#([n]) \to (\Coker\phi)([n]) \to 0.
$$
The transtition maps $\iota_f$ for $\Ker\phi$ are obtained by
restriction from those for $M_\#$; for $\Coker\phi$, one uses the
fact that the functors $f_!$ are right-exact.
\endproof

\begin{defn}
A {\em cyclic bimodule $M$} over a unital associative algebra $A$ is
a cocartesian section $M_\# \in \Sec_{cart}(A\bimod_\#)$. A {\em
complex of cyclic bimodules $M_\idot$} over $A$ is an object in the
derived category $\D(\Sec(A\bimod_\#))$ whose homology objects are
cocartsian.
\end{defn}

Complexes of cyclic bimodules obviously form a full triangulated
subcategory in $\D(\Sec(A\bimod_\#))$; consistent notation for this
category would be $\D_{cart}(\Sec(A\bimod_\#))$, but for simplicity
we will denote it $\D\Lambda(A\bimod)$. We have to define complexes
separately for the following reasons:
\begin{enumerate}
\item The category $\Sec_{cart}(A\bimod_\#) \subset
  \Sec(A\bimod_\#)$ need not be abelian -- since the transition
  functors $f_!$ are only right-exact, the condition of being
  cocartesian need not be preserved when passing to kernels.
\item Even if $\Sec_{cart}(A\bimod_\#)$ is abelian, its derived
  category might be much smaller than $\D\Lambda(A\bimod)$.
\end{enumerate}

\begin{example}\label{const.exa}
An extreme example of \thetag{ii} is the case $A = k$: in this case
$\Sec(A\bimod_\#)$ is just the category of cyclic vector spaces,
$\Fun(\Lambda,k)$, and $E \in \Fun(\Lambda,k)$ is cocartesian if and
only if $E(f)$ is invertible for any map $f:[n'] \to [n]$. One
deduces easily that $E$ must be a constant functor, so that
$\Sec_{cart}(k\bimod_\#) = k\Vect$. Then $\D\Lambda(k\bimod)$ is the
full subcategory $\D_{const}(\Lambda,k) \subset \D(\Lambda,k)$ of
complexes whose homology is constant. If we were to consider
$\Delta^{opp}$ instead of $\Lambda$, we would have
$\D_{const}(\Delta^{opp},k) \cong \D(k\Vect)$ -- since
$H^\hdot(\Delta^{opp},k) = k$, the embedding $\D(k\Vect) \to
\D(\Delta^{opp},k)$ is fully faithful, and
$\D_{const}(\Delta^{opp},k)$ is its essential image. However,
$H^\hdot(\Lambda,k)$ is $k[u]$, not $k$. Therefore there are maps
between constant functors in $\D(\Lambda,k)$ which do not come from
maps in $\D(k\Vect)$, and the cones of these maps give objects in
$\D_{const}(\Lambda,k)$ which do not come from $\D(k\Vect)$.
\end{example}

This phenomenon is quite common in homological algebra -- examples
are, for instance, the triangulated category of complexes of \'etale
sheaves with constructible homology, the category of complex of
$\D$-modules with holonomic homology, or the so-called ``equivariant
derived category'' of sheaves on a topological space $X$ acted upon
by a topological group $G$ (which is not in fact the derived
category of anything useful). The upshot is that it is the
triangulated category $\D\Lambda(A\bimod)$ which should be treated
as the basic object, wherever categories are discussed.

\begin{remark}\label{const.rem}
We note one interesting property of the category
$\D_{const}(\Lambda,k)$. Fix an integer $n \geq 1$, and consider the
full subcategory $\Lambda_{\leq n} \subset \Lambda$ of objects $[n']
\in \Lambda$ with $n' \leq n$. Then one can show that
$H^\hdot(\Lambda_{\leq n},k) = k[u]/u^n$, so that we have a natural
exact triangle
\begin{equation}\label{conn.2}
\begin{CD}
H_\idot(\Lambda_{\leq n},E^\hdot) @>>> HC_\idot(E^\hdot) @>{u^n}>>
HC_{\idot+2n}(E) @>>>,
\end{CD}
\end{equation}
for every $E^\hdot \in \D_{const}(\Lambda,k)$. We note that for any
$E^\hdot \in \D(\Lambda,k)$, \eqref{connes} extends to a spectral
sequence
\begin{equation}\label{conn.sp}
HH_\idot(E^\hdot)[u^{-1}] \Rightarrow HC_\idot(E),
\end{equation}
where the expression on the left-hand side reads as ``polynomials in
one formal variable $u^{-1}$ of homological degree $2$ with
coefficients in $HH_\idot(E^\hdot)$''. Then \eqref{conn.2} shows
that for $E^\hdot \in \D_{const}(\Lambda,k)$, the first $n$
differentials in \eqref{conn.sp} depend only on the restriction of
$E^\hdot$ to $\Lambda_{\leq (n+1)} \subset \Lambda$. This is useful
because in practice, one is often interested only in the first
differential in the spectral sequence.
\end{remark}

As in Lemma~\ref{cycl.str}, a cyclic $A$-bimodule $M_\#$ essentially
consists of an $A$-bimodule $M = M_\#([1])$ equipped with an extra
structure. Explicitly, this structure is a map $\tau:A \otimes_k M
\to M \otimes_k A$ which respects the $A^{\otimes 2}$-bimodule
structure on both sides, and satisfies the condition $\tau_{31}
\circ \tau_{12} \circ \tau_{23} = \id$, as in Lemma~\ref{cycl.str}.

Another way to view this structure is the following. One checks
easily that for any cyclic $A$-bimodule $M_\#$, the restriction
$j^*M_\# \in \Fun(\Delta^{opp},k)$ is canonically isomorphic to the
simplicial $k$-vector space $M^\Delta_\#$ associated to the
underlying $A$-bimodule $M$ as in \eqref{M.Delta}. By adjunction, we
have a natural map
$$
\tau_\#:j_!M^\Delta_\# \to M_\#.
$$
Then $j_!M^\Delta_\#$ in this formula only depends on $M \in
A\bimod$, and all the structure maps which turn $M$ into the cyclic
bimodule $M_\#$ are collected in the map $\tau_\#$.

\medskip

We can now define cyclic homology with coefficients. The definition
is rather tautological. We note that for any cyclic $A$-bimodule
$M_\#$ -- or in fact, for any $M_\# \in \Sec(A\bimod_\#)$ -- we can
treat $M_\#$ as a cyclic vector space by forgetting the bimodule
structure on its components $M_n$.

\begin{defn}\label{cycl.def}
The {\em cyclic homology $HC_\idot(A,M_\#)$ with coefficients} in a
cyclic $A$-bimodule $M$ is equal to $H_\idot(\Lambda,M_\#)$.
\end{defn}

Of course, \eqref{connes}, being valid for any cyclic $k$-vector
space, also applies to $HC_\idot(A,M_\#)$, so that we automatically
get the whole package -- the Connes' exact sequence, the periodicity
endomorphism, and the periodic cyclic homology $HP_\idot(A,M)$. By
Lemma~\ref{hoch}, $HH_\idot(M_\#)$ coincides with $HH_\idot(A,M)$ as
defined in \eqref{hoch.coeff}.

\section{Gauss-Manin connection.}\label{getz}

To illustate the usefulness of the notion of a cyclic bimodule, let
us study the behavior of cyclic homology under deformations.

There are two types of deformation theory objects that one can study
for an associative algebra $A$. The first is the notion of a {\em
square-zero extension} of the algebra $A$ by a $A$-bimodule
$M$. This is an associative algebra $\wt{A}$ which fits into a short
exact sequence
$$
\begin{CD}
0 @>>> M @>{i}>> \wt{A} @>{p}>> A @>>> 0,
\end{CD}
$$
where $p$ is an algebra map, and $i$ is an $\wt{A}$-bimodule map,
under the $\wt{A}$-bimodule structure on $M$ induced from the given
$A$-bimodule structure by means of the map $p$. In other words, the
multiplication on the ideal $\Ker p \subset \wt{A}$ is trivial, so that
the $\wt{A}$-bimodule structure on $\Ker p$ is induced by an
$A$-bimodule structure, and $i$ identifies the $A$-bimodule $\Ker p$
with $M$. Square-zero extensions are classified up to an isomorphism
by elements in the second Hochschild cohomology group $HH^2(A,M)$,
defined as
$$
HH^\hdot(A,M) = \Ext^\hdot_{A^{opp} \otimes A}(A,M).
$$
In this setting, we can consider the cyclic homology of the algebra
$\wt{A}$ and compare with the cyclic homology of $A$.
Th.~Goodwillie's theorem \cite{go} claims that if the base field $k$
has characteristic $0$, the natural map
$$
HP_\idot(\wt{A}) \to HP_\idot(A)
$$
is an isomorphism, and there is also some information on the
behaviour of $HC_\idot(A)$.

A second type of deformation theory data includes a commutative
$k$-algebra $R$ with a maximal ideal $\m \subset R$. A {\em
deformation} $A_R$ of the algebra $A$ over $R$ is a flat associative
unital algebra $A_R$ over $R$ equipped with an isomorphism $A_R/\m
\cong A$. In this case, one can form the {\em relative} cyclic
$R$-module $A_{R\#}$ by taking the tensor products over $R$; thus we
have relative homology $HH_\idot(A_R/R)$, $HC_\idot(A_R/R)$,
$HP_\idot(A_R/R)$. The fundamental fact discovered by E. Getzler
\cite{getz} is that we have an analog of the Gauss-Manin connection:
if $\Spec R$ is smooth, the $R$-module $HP_i(A_R/R)$ carries a
canonical flat connection for every $i$.

Consider now the case when $R$ is not smooth but, on the contrary,
local Artin. Moreover, assume that $\m^2=0$, so that $R$ is itself a
(commutative) square-zero extension of $k$. Then a deformation $A_R$
of $A$ over $R$ is also a square-zero extension of $A$, by the
bimodule $A \otimes \m$ ($\m$ here is taken as a $k$-vector
space). But this square-zero extension is special -- for a general
square-zero extension $\wt{A}$ of $A$ by some $M \in A\bimod$, there
does not exist any analog of the relative cyclic $R$-module $A_{R\#}
\in \Fun(\Lambda,R)$.

\medskip

We observe the following: the data needed to define such an analog
is precisely a cyclic bimodule structure on the bimodule $M$.

\medskip

Namely, assume given a square-zero extension $\wt{A}$ of the algebra
$A$ by some $A$-bimodule $M$, and consider the cyclic $k$-vector
space $\wt{A}_\# \in \Fun(\Lambda,k)$. Let us equip $\wt{A}$ with a
descreasing two-step filtration $F^\hdot$ by setting $F^1\wt{A} =
M$. Then this induces a decreasing filtration $F^\hdot$ on tensor powers
$\wt{A}^{\otimes n}$. Since $\wt{A}$ is square-zero, $F^\hdot$ is
compatible with the multiplication maps; therefore we also have a
filtration $F^\hdot$ on $\wt{A}_\#$. Consider the quotient
$$
\overline{A_\#} = \wt{A}_\#/F^2\wt{A}_\#.
$$
One checks easily that $\gr^0_F\wt{A}_\# \cong A_\#$ and
$\gr^1_F\wt{A}_\# \cong j_!M^\Delta_\#$ in a canonical way, so that
$\overline{A_\#}$ fits into a canonical short exact sequence
\begin{equation}\label{A.bar}
\begin{CD}
0 @>>> j_!M^\Delta_\# @>>> \overline{A_\#} @>>> A_\# @>>> 0
\end{CD}
\end{equation}
of cyclic $k$-vector spaces.

Now assume in addition that $M$ is equipped with a structure of a
cyclic $A$-bimodule $M_\#$, so that $M^\Delta_\# \cong
j^*M_\#$, and we have the structure map $\tau_\#:j_!M^\Delta_\# \to
M_\#$. Then we can compose the extension \eqref{A.bar} with the map
$\tau_\#$, to obtain a commutative diagram
\begin{equation}\label{A.hat}
\begin{CD}
0 @>>> j_!M^\Delta_\# @>>> \overline{A_\#} @>>> A_\# @>>> 0\\
@. @V{\tau_\#}VV @VVV @|\\
0 @>>> M_\# @>>> \wh{A_\#} @>>> A_\# @>>> 0
\end{CD}
\end{equation}
of short exact sequences in $\Fun(\Lambda,k)$, with cartesian left
square. It is easy to check that when $\wt{A} = A_R$ for some
square-zero $R$, so that $M = A \otimes \m$, and we take the
cyclic $A$-bimodule structure on $M$ induced by the tautological
structure on $A$, then $\wh{A_\#}$ coincides precisely with the
relative cyclic object $A_{R\#}$ (which we consider as a $k$-vector
space, forgetting the $R$-module structure).

We believe that this is the proper generality for the Getzler
connection; in this setting, the main result reads as follows.

\begin{prop}\label{spl}
Assume given a square-zero extension $\wt{A}$ of an associative
algebra $A$ by an $A$-bimodule $M$, and assume that $M$ is equipped
with a structure of a cyclic $A$-bimodule. Then the long exact
sequence
$$
\begin{CD}
HP_\idot(A,M) @>>> HP_\idot(\wh{A_\#}) @>>> HP_\idot(A) @>>>
\end{CD}
$$
of periodic cyclic homology induced by the second row in
\eqref{A.hat} admits a canonical splitting $HP_\idot(A) \to
HP_\idot(\wh{A_\#})$.
\end{prop}

\proof{} By definition, we have two natural maps
\begin{equation}\label{comp}
\begin{aligned}
HP_\idot(\overline{A_\#}) &\to HP_\idot(A_\#) = HP_\idot(A),\\
HP_\idot(\overline{A_\#}) &\to HP_\idot(\wh{A_\#}),
\end{aligned}
\end{equation}
and the cone of the first map is isomorphic to
$HP_\idot(j_!M^\Delta_\#)$. Since $j_!$ is exact, we have
$HC_\idot(j_!M^\Delta_\#) \cong HH_\idot(M_\#)$, and the periodicity
map $u:HC_\idot(j_!M^\Delta_\#) \to HC_{\idot-2}(j_!M^\Delta_\#)$ is
equal to $0$, so that $HP_\idot(j_!M^\Delta_\#) = 0$. Thus the first
map in \eqref{comp} is an isomorphism, and the second map is then
the required splitting.
\endproof

\begin{corr}
Assume given a commutative $k$-algebra $R$ with a maximal ideal $\m
\subset R$, and a deformation $A_R$ of the algebra $A$ over
$R$. Then if $\Spec R$ is smooth, the $R$-modules $HP_\idot(A_R/R)$
carry a natural connection.
\end{corr}

\proof[Sketch of a proof.] Consider the $R \otimes R$-algebras $A_R
\otimes R$ and $R \otimes A_R$, and their restrictions to the first
infinitesemal neighborhood of the diagonal in $\Spec (R \otimes R) =
\Spec R \times \Spec R$. Then Proposition~\ref{spl}, suitably
generalized, shows that $HP_\idot(-)$ of these two restrictions are
canonically isomorphic. It is well-known that giving such an
isomorphism is equivalent to giving a connection on
$HP_\idot(A_R/R)$.
\endproof

We note that we do not claim that the connection is {\em flat}. It
certainly is, at least in characteristic $0$; but our present method
does not allow one to go beyond square-zero extensions. Thus we
cannot analyse the second infinitesemal neighborhood of the diagonal
in $\Spec (R \otimes R)$, and we cannot prove
flatness.

Unfortunately, at present, we do not understand what is the proper
cyclic bimodule context for higher-level infinitesemal
extensions. Of course, if one is only interested in an
$R$-deformation $\wt{A} = A_R$ over an Artin local base $R$, not in
its cyclic bimodule generalizations, one can use Goodwillie's
Theorem: using the full cyclic object $\wt{A}_\#$ instead of its
quotient $\overline{A_\#}$ in Proposition~\ref{spl} immediately
gives a splitting $HP_\idot(A) \to HP_\idot(A_R/R)$ of the
augmentation map $HP_\idot(A_R/R) \to HP_\idot(A)$, and this extends
by $R$-linearity to an isomorphism $HP_\idot(A_R/R) \cong
HP_\idot(A) \otimes R$. However, this is not quite satisfactory from
the conceptual point of view, and it does not work in positive
characteristic (where Goodwillie's Theorem is simply not true). If
$\cchar k \neq 2$, the latter can be cured by using
$\wt{A}_\#/F^3\wt{A}_\#$, but the former remains. We plan to return
to this elsewhere.

\section{Categorical approach.}\label{cat}

Let us now try to define cyclic homology in a more general setting
-- we will attempt to replace $A\bimod$ with an arbitrary
associative unital $k$-linear tensor category $\C$ with a unit
object $\I \in \C$. We do not assume that $\C$ is symmetric in any
way. However, we will assume that the tensor product $- \otimes -$
is right-exact in each variable, and we will need to impose
additional technical assumptions later on.

The first thing to do is to try to define Hochschild homology; so,
let us look more closely at \eqref{hh.def}. The formula in the
right-hand side looks symmetric, but this is an optical illusion --
the two copies of $A$ are completely different objects: one is a
left module over $A^{opp} \otimes A$, and the other is a right
module ($A$ just happens to have both structures at the same
time). It is better to separate them and introduce the functor
$$
\tr:A\bimod \to k\Vect
$$
by $\tr(M) = M \otimes_{A^{opp} \otimes A} A$ -- or, equivalently,
by
\begin{equation}\label{tr.A}
\tr(M) = M/\{ am-ma \mid a \in A, m \in M \}.
\end{equation}
Then $\tr$ is a right-exact functor, and we have $HH_\idot(A,M) =
L^\hdot\tr(M)$. 

We want to emphasize that the functor $\tr$ can not be recovered
from the tensor structure on $A\bimod$ -- this really is an extra
piece of data. For a general tensor category $\C$, it does not exist
a priori; we have to impose it as an additional structure.

Let us axiomatize the situation. First, forget for the moment about
the $k$-linear and abelian structure on $\C$ -- let us treat it
simply as a monoidal category. Assume given some other category $\B$
and a functor $T:\C \to \B$.

\begin{defn}\label{trace.defn}
The functor $T:\C \to \B$ is a {\em trace functor} if it is extended
to a functor $\C_\# \to \B$ which sends any cocartesian map $f:M \to
M'$ in $\C_\#$ to an invertible map.
\end{defn}

Another way to say the same thing is the following: the categories
$\Fun(\C^n,\B)$ of functors from $\C^n$ to $\B$ form a fibered
category over $\Lambda$, and a trace functor is a cartesian section
of this fibration. Explicitly, a trace functor is defined by $T:\C
\to \B$ and a collection of isomorphisms
$$
T(M \otimes M') \to T(M' \otimes M)
$$
for any $M,M' \in \C$ which are functorial in $M$ and $M'$ and
satisfy some compatibility conditions analogous to those in
Lemma~\ref{cycl.str}; we leave it to the reader to write down these
conditions precisely. Thus $T$ has a trace-like property with
respect to the product in $\C$, and this motivates our terminology.

Recall now that $\C$ is a $k$-linear abelian category. To define
Hochschild homology, we have to assume that it is equipped with a
right-exact trace functor $\tr:\C \to k\Vect$; then for any $M \in
\C$, we set
\begin{equation}\label{hh.def.gen}
HH_\idot(M) = L^\hdot\tr(M).
\end{equation}

\begin{lemma}
The functor $\tr:A\bimod \to k\Vect$ canonically extends to a
right-exact trace functor in the sense of
Definition~\ref{trace.defn}.
\end{lemma}

\proof{} For any object $\langle [n],M_n \rangle \in A\bimod_\#$,
$[n] \in \Lambda$, $M_n \in A^{\otimes n}\bimod$, let
$$
\tr(\langle [n],M_n \rangle) = M_n/\{ a_{v'}m - ma_v \mid v \in
V([n]), m \in M_n, a \in A \},
$$
where $a_v = 1 \otimes 1 \otimes \dots \otimes a \otimes \dots
\otimes 1 \in A^{\otimes V([n])}$ has $a$ in the multiple
corresponding to $v \in V([n])$, and $v' \in V([n])$ is the next
marked point after $v$ counting clockwise. The compatibility with
maps in the category $A\bimod_\#$ is obvious.
\endproof

We note that here, in the case $\C=A\bimod$, the category
$A\bimod_\#$ is actually larger than what we would have had purely
from the monoidal structure on $\C$: $M_n$ is allowed to be an
arbitrary $A^{\otimes n}$-bimodule, not a collection of $n$
$A$-bimodules. To do the same for general $k$-linear $\C$, we need
to replace $A^{\otimes n}\bimod$ with some version of the tensor
product $\C^{\otimes n}$. Here we have a difficulty: for various
technical reasons, it is not clear how to define tensors products
for sufficiently general abelian categories.

One way around it is the following. For any (small) $k$-linear
abelian category $\B$, a $k$-linear functor $\B^{opp} \to k\Vect$ is
left-exact if and only if it is a sheaf for for the canonical
Grothendieck topology on $\B$ (\cite[5, \S10]{BD}); the category
$\Shv(\B)$ of such functors is abelian and $k$-linear, and $\B$
itself is naturally embedded into $\Shv(\B)$ by Yoneda. The
embedding is a fully faithfull exact functor. Every functor in
$\Shv(\B)$ is in fact a direct limit of representable functors, so
that $\Shv(\B)$ is an inductive completion of the abelian category
$\B$. Now, if are given two (small) $k$-linear abelian categories
$\B_1$, $\B_2$, then their product $\B_1 \times \B_2$ is no longer
abelian. However, we still have the abelian category $\Shv(\B_1
\times \B_2)$ of bilinear functors $\B_1^{opp} \times \B_2^{opp} \to
k\Vect$ which are left-exact in each variable, and the same goes for
polylinear functors.

Moreover, for any right-exact functor $F:\B_1 \to \B_2$ between
small abelian categories, we have the restriction functor
$F^*:\Shv(\B_2) \to \Shv(\B_1)$, which is left-exact, and its
left-adjoint $F_!:\Shv(\B_1) \to \Shv(\B_2)$, which is
right-exact. The functor $F_!$ is an extension of the functor $F$:
on Yoneda images $\B_i \subset \Shv(\B_i)$, we have $F_! = F$. And,
again, the same works for polylinear functors.

In particular, given our $k$-linear abelian tensor category $\C$, we
can form the category $\Shv(\C)_\#$ of pairs $\langle E,[n]
\rangle$, $[n] \in \Lambda$, $E \in \Shv(\C^n)$, with a map from
$\langle E',[n'] \rangle$ to $\langle E,[n] \rangle$ given by a pair
of a map $f:[n'] \to [n]$ and either a map $E' \to (f_!)^*E$, or map
$(f_!)_!E' \to E$ -- this is equivalent by adjunction. Then
$\Shv(\C)_\#$ is bifibered category over $\Lambda$ in the sense of
\cite{SGA}.

The category of sections $\Lambda \to \Shv(\C)_\#$ of this
bifibration can also be described as the full subcategory
$\Shv(\C_\#) \subset \Fun(\C^{opp}_\#,k)$ spanned by those functors
$E_\#:\C^{opp}_\# \to k\Vect$ whose restriction to $(\C^{opp})^n
\subset \C^{opp}_\#$ is a sheaf -- that is, an object in $\Shv(\C^n)
\subset \Fun((\C^{opp})^n,k)$. Since the transition functors
$(f_!)_!$ are right-exact, $\Shv(\C_\#)$ is an abelian category
(this is proved in exactly the same way as Lemma~\ref{sec.ab}).

We denote by $\Shv_{cart}(\C_\#) \subset \Shv(\C_\#)$ the full
subcategory of sections $E:\Lambda \to \Shv(\C)_\#$ which are
cocartesian, and moreover, are such that $E([1]) \in \Shv(\C)$
actually lies in the Yoneda image $\C \subset \Shv(\C)$. We also
denote by $\D\Lambda(\C) \subset \D(\Shv(\C_\#))$ the full
triangulated subcategory of complexes $E^\hdot_\# \in
\D(\Shv(\C_\#))$ with homology in $\Shv_{cart}(\C_\#)$.

\medskip

If $\C$ is the category of $A$-bimodules for some algebra $A$ -- or
better yet, of $A$-bimodules of cardinality not more than that of $A
\times \N$, so that $\C$ is small -- then $\Shv(\C)$ is equivalent
to $A\bimod$ (one shows easily that every sheaf $E \in \Shv(\C)$ is
completely determined by its value at $A^{opp} \otimes A \in
\C$). In this case, $\D\Lambda(\C)$ is our old category
$\D\Lambda(A\bimod)$.

\medskip

Now, we assume that $\C$ is equipped with a right-exact trace
functor $\tr:\C \to k\Vect$, we would like to define cyclic homology
$HC_\idot(M_\idot)$ for any $M_\idot \in \D\Lambda(\C)$, and we
immediately notice a problem: for a general $\C$, we do not have a
forgetful functor to vector spaces. However, it turns out that the
forgetful functor {\em is not needed} for the definition -- it can
be replaced with the trace functor $\tr$.

We proceed as follows. By definition, $\tr$ is extended to a functor
$\C_\# \to k\Vect$; we extend it canonically to a functor
$\Shv(\C)_\# \to k\Vect$, and consider the product
$$
\tr \times \tau: \Shv(\C)_\# \to k\Vect \times \Lambda,
$$
where $\tau:\Shv(\C)_\# \to \Lambda$ is the projection.  This is a
functor compatible with the projections to $\Lambda$, and therefore,
it induces a functor of the categories of sections. The category of
sections of the projection $k\Vect \times \Lambda \to \Lambda$ is
tautologically the same as $\Fun(\Lambda,k\Vect)$, so that we have a
functor
$$
\tr_\#:\Shv(\C_\#) \to \Fun(\Lambda,k).
$$
One checks easily that this functor is right-exact. 

\begin{defn}\label{cycl.def.gen}
For any $M_\# \in \Sec(\C_\#)$, its cyclic homology $HC_\idot(M_\#)$
is defined by
$$
HC_\idot(M_\#) \defeq HC_\idot(L^\hdot\tr_\#(M_\#)) =
H_\idot(\Lambda,L^\hdot\tr_\#(M_\#)).
$$
\end{defn}

\begin{defn}\label{clean}
The pair $\langle \C,\tr \rangle$ is called {\em homologically
clean} if for any $n$, the category $\Shv(\C^n)$ has enough objects
$E$ such that
\begin{enumerate}
\item $E$ is acyclic both for functors $(f_!)_!:\Shv(\C^n) \to
\Shv(\C^{n'})$, for any $f:[n] \to [n']$, and for the trace functor
$\tr:\Shv(\C^n) \to k\Vect$, and
\item for any $f:[n] \to [n']$, $(f_!)_!E \in \Shv(\C^{n'})$ is
  acyclic for $\tr:\Shv(\C^{n'}) \to k\Vect$.
\end{enumerate}
\end{defn}

\begin{example}\label{cln.exa}
Assume that the category $\C$ has enough projectives, and moreover,
$P_1 \otimes P_2$ is projective for any projective $P_1,P_2 \in \C$
(this is satisfied, for instance, for $\C=A\bimod$). Then the pair
$\langle \C,\tr \rangle$ is homologically clean, for any trace
functor $\tr$. Indeed, $\Shv(\C^n)$ then also has enough
projectives, say sums of objects
\begin{equation}\label{box}
P = P_1 \boxtimes P_2 \boxtimes \dots \boxtimes P_n \in \Shv(\C^n)
\end{equation}
for projective $P_1,\dots,P_n \in \C \subset \Shv(\C)$, and these
projectives automatically satisfy the condition \thetag{i}. To check
\thetag{ii}, one decomposes $f:[n] \to [n']$ into a surjection
$p:[n] \to [n'']$ and an injection $i:[n''] \to [n']$. Since the
tensor product of projective objects is projective, $(p_!)_!(P) \in
\Shv(\C^{n''})$ is also an object of the type \eqref{box}, so we may
as well assume that $f$ is injective. Then one can find a
left-inverse map $f':[n'] \to [n]$, $f' \circ f = \id$; since
$P'=(f_!)_!(P)$ is obviously acyclic for $(f'_!)_!$, and
$(f'_!)_!(P') = ((f' \circ f)_!)_!(P) = P$ is acyclic for $\tr$,
$P'$ itself is acyclic for $\tr = \tr \circ (f'_!)_!$.
\end{example}

\begin{lemma}\label{b.ch}
Assume that $\langle \C,\tr \rangle$ is homologically clean. Then
for any object $[n] \in \Lambda$ and any $M_\# \in \Shv(\C_\#)$, we
have
\begin{equation}\label{equa}
L^\hdot\tr_\#(M_\#)([n]) \cong L^\hdot\tr (M_\#([n])).
\end{equation}
For any $M_\#^\hdot \in \D\Lambda(\C)$, we have $L^\hdot\tr_\#(M_\#)
\in \D_{const}(\Lambda,k) \subset \D(\Lambda,k)$.
\end{lemma}

\proof{} The natural restriction functor $\Shv(\C_\#) \to
\Shv(\C^n)$, $M_\# \mapsto M_\#([m])$ has a left-adjoint functor
$I_{n!}:\Shv(\C^n) \to \Shv(\C_\#)$; explicitly, it is given by
\begin{equation}\label{ind}
I_{n!}(E)([n']) = \bigoplus_{f:[n] \to [n']}(f_!)_!(E).
\end{equation}
Let us say that an object $E \in \Shv(\C^n)$ is admissible if it
satisfies the conditions \thetag{i}, \thetag{ii} of
Definition~\ref{clean}. By assumption, $\Shv(\C^n)$ has enough
admissible objects for any $n$. Then $\Shv(\C_\#)$ has enough
objects of the form $I_{n!}E$, $[n] \in \Lambda$, $E \in \Shv(\C^n)$
admissible, and to prove the first claim, it suffices to consider
$M_\#=I_{n!}E$ of this form. In degree $0$, \eqref{equa} is the
definition of the functor $\tr_\#$, and the higher degree terms in
the right-hand side vanish by Definition~\ref{clean}~\thetag{ii}.
Therefore it suffices to prove that $M_\#=I_{n!}E$ is acyclic for
the functor $\tr_\#$. This is obvious: applying $\tr_\#$ to any
short exact sequence
$$
\begin{CD}
0 @>>> M'_\# @>>> M''_\# @>>> M_\# @>>> 0
\end{CD}
$$
in $\Shv(\C_\#)$, we see that, since $M_\#([n'])$ is acyclic for any
$[n'] \in \Lambda$, the sequence
$$
\begin{CD}
0 @>>> \tr M'_\#([n']) @>>> \tr M''_\# ([n']) @>>> \tr M_\# ([n']) @>>> 0
\end{CD}
$$
is exact; this means that
$$
\begin{CD}
0 @>>> \tr M'_\# @>>> \tr M''_\# @>>> \tr M_\# @>>> 0
\end{CD}
$$
is an exact sequence in $\Fun(\Lambda,k)$, and this means that
$M_\#$ is indeed acyclic for $\tr_\#$.

With the first claim proved, the second amounts to showing that the
natural map
$$
L^\hdot\tr \circ L^\hdot(f_!)_!(E) \to L^\hdot\tr(E)
$$
is a quasiismorphism for any $f:[n] \to [n']$ and any $E \in
\Shv(\C^n)$. It suffices to prove it for admissible $M$; then the
higher derived functors vanish, and the isomorphism $\tr \circ
(f_!)_! \cong \tr$ is Definition~\ref{trace.defn}.
\endproof

\begin{lemma}
In the assumptions of Lemma~\ref{b.ch}, for any complex $M^\hdot_\#
\in \D\Lambda(\C)$ with the first component $M^\hdot =
M^\hdot_\#([1])$ we have
$$
HH_\idot(M^\hdot) \cong HH_\idot(L^\hdot\tr_\#(M^\hdot_\#)).
$$
\end{lemma}

\proof{} By Lemma~\ref{b.ch}, the left-hand side,
$HH_\idot(M^\hdot)$, is canonically isomorphic to the complex
$L^\hdot\tr_\#(M^\hdot_\#) \in \D(\Lambda,k)$ evaluated at $[1] \in
\Lambda$, and moreover, $L^\hdot\tr_\#(M^\hdot_\#)$ lies in the
subcategory $\D_{const}(\Lambda,k) \subset \D(\Lambda,k)$. It
remains to apply the general fact: for any $E^\hdot \in
\D_{const}(\Lambda,k)$, we have a natural isomorphism
$HH_\idot(E^\hdot) \cong E^\hdot([1])$. Indeed, by definition we
have
$$
HH_\idot(E^\hdot) = H_\idot(\Delta^{opp},j^*E^\hdot),
$$
and $j^*E^\hdot$ lies in the category $\D_{const}(\Delta^{opp},k)$
which is equivalent to $\D(k\Vect)$ (see Example~\ref{const.exa}, and
also Remark~\ref{const.rem}: the isomorphism we constructed here is
a special case of \eqref{conn.2} for $n=1$).
\endproof

The Lemma shows that if the pair $\langle \C,\tr \rangle$ is
homologically clean, Definition~\ref{cycl.def.gen} is consistent
with \eqref{hh.def.gen}, and we get the whole periodicity package of
\eqref{connes} -- the periodicity map $u$, the Connes' exact
sequence
$$
\begin{CD}
HH_\idot(M^\hdot) @>>> HC_\idot(M^\hdot) @>{u}>> HC_{\idot-2}(M^\hdot) @>>>,
\end{CD}
$$
and the periodic cyclic homology $HP_\idot(M^\hdot)$.

In general, objects in $\D\Lambda(\C)$ may be hard to construct, but
we always have at least one -- the identity section $\I_\#:\Lambda
\to \Shv(\C)_\#$, given by
$$
\I_\#([n]) = \I^{\boxtimes n} \in \C^{\otimes n},
$$
where $\I \in \C$ is the unit object. Thus we can define cyclic
homology of a tensor category equipped with a trace functor.

\begin{defn}\label{cycl.cat}
For any $k$-linear abelian unital tensor category $\C$ equipped with
a trace functor $\tr:\C \to k\Vect$, its Hochschild and cyclic
homology is given by
$$
HH_\idot(\C,\tr) \defeq HH_\idot(\I), \qquad
HC_\idot(\C,\tr) \defeq HC_\idot(\I_\#),
$$
where $\I \in \C$ is the unit object, and $\I_\# \in \D\Lambda(\C)$
is the identity section.
\end{defn}

We now have to check that in the case $\C = A\bimod$,
Definition~\ref{cycl.def.gen} is compatible with our earlier
Definition~\ref{cycl.def} -- in other words, that the cyclic
homology computed by means of the forgetfull functor is the same as
the cyclic homology computed by means of the trace.  This is not at
all trivial. Indeed, if for instance $M_\# \in \Shv(\C_\#)$ is
cocartesian, then, while $L^\hdot\tr^\#M_\#$ lies in the subcategory
$\D_{const}(\Lambda,k) \subset \D(\Lambda,k)$, the same is certainly
not true for the object $M_\# \in \Fun(\Lambda,k)$ obtained by
forgetting the bimodule structure on $M_n$.

Thus these two objects are different. However, they do become equal
after taking cyclic (or Hochschild, or periodic cyclic)
homology. Namely, for any $M_\# \in \Sec(A\bimod_\#)$ we have a
natural map
\begin{equation}\label{natu}
M_\# \to L^\hdot\tr^\#M_\#
\end{equation}
in the derived category $\D(\Lambda,k)$, and we have the following
result.

\begin{prop}\label{main}
For every $M_\# \in \Sec(A\bimod_\#)$, the natural map \eqref{natu}
induces isomorphisms
\begin{align*}
HH_\idot(M_\#) &\cong HH_\idot(L^\hdot\tr M_\#),\\
HC_\idot(M_\#) &\cong HC_\idot(L^\hdot\tr M_\#),\\
HP_\idot(M_\#) &\cong HP_\idot(L^\hdot\tr M_\#).
\end{align*}
\end{prop}

\proof{} By \eqref{connes}, it suffices to consider $HC_\idot(-)$;
as in the proof of Lemma~\ref{b.ch}, it suffices to consider $M_\# =
I_{n!}E$ given in \eqref{ind}, with $E$ being the free bimodule
$$
E = (A^{opp} \otimes A)^{\otimes n} \in \Shv(\C^n) = A^{\otimes
  n}\bimod
$$
for some fixed $n$. Explicitly, we have
\begin{equation}\label{ind.2}
I_{n!}E([n']) = \bigoplus_{f:[n] \to [n']}\bigotimes_{v' \in V([n'])}
A^{opp} \otimes A^{\otimes f^{-1}(v')}
\end{equation}
for any $[n'] \in \Lambda$. Then $L^p\tr_\#I_{n!}E = 0$ for $p \geq 1$,
and one checks easily that
$$
\tr_\# I_{n!}E = i_{n!}\tr E = i_{n!}A^{\otimes n} \in \Fun(\Lambda,k),
$$
where $i_n:\ppt \to \Lambda$ is the embedding of the object $[n] \in
\Lambda$ ($\ppt$ is the category with one object and one
morphism). Therefore
$$
HC_0(L^\hdot\tr_\#I_{n!}E) = H_\idot(\Lambda,i_{n!}A^{\otimes n}) =
A^{\otimes n},
$$
and $HC_p(L^\hdot\tr_\#i_{n!}E) = 0$ for $p \geq 1$. We have to compare
it with $HC_\idot(i_{n!}E)$.

To do this, consider the category $\Lambda_{[n]}$ of objects $[n']
\in \Lambda$ equipped with a map $[n] \to [n']$, and let
$j_n:\Lambda_{[n]} \to \Lambda$ be the forgetful functor. Then $j_n$
is obviously a discrete cofibration. Comparing \eqref{discr} and
\eqref{ind.2}, we see that
$$
I_{n!}E = j_{n!}E_\#^{[n]}
$$
for some $E_\#^{[n]} \in \Fun(\Lambda_{[n]})$. Moreover, fix once
and for all a map $[1] \to [n]$. Then we see that the discrete
cofibration $j_n:\Lambda_{[n]} \to \Lambda$ factors through the
discrete cofibration $j:\Lambda_{[1]} = \Delta^{opp} \to \Lambda$ by
means of a discrete cobifbration $\gamma_n:\Lambda_{[n]} \to
\Lambda_{[1]}$, and we observe that
$$
E^{[n]}_\#([n']) = (A^{opp})^{\otimes n'} \otimes A^{\otimes n}
$$
only depends on $\gamma_n([n']) \in \Delta^{opp}$. More precisely,
we have $E^{[n]}_\# = \gamma^*_nE_n^{\Delta}$, where $E_n^{\Delta}
\in \Fun(\Delta^{opp},k)$ is as in \eqref{M.Delta}, and $E_n$ is the
free $A$-bimodule
$$
E_n = A^{opp} \otimes A^{\otimes (n-1)} \otimes A.
$$
The conclusion: we have
$$
HC_\idot(I_{n!}E) = H_\idot(\Lambda_{[n]},E^{[n]}_\#)
= H_\idot(\Delta^{opp},\gamma_{n!}\gamma^*_nE_n^{\Delta}) =
H_\idot(\Delta^{opp},E_n^{\Delta} \otimes \gamma_{n!}k),
$$
where we have used the projection formula \eqref{projj} in the
right-hand side. The homology of the category $\Delta^{opp}$ can be
computed by the standard complex; then by the K\"unneth formula, the
right-hand side is isomorphic to
$$
H_\idot(\Delta^{opp},E_n^{\Delta}) \otimes
H_\idot(\Delta^{opp},\gamma_{n!}k) 
\cong H_\idot(\Delta^{opp},E_n^{\Delta}) \otimes
H_\idot(\Lambda_{[n]},k).
$$
By Lemma~\ref{hoch},
$$
H_\idot(\Delta^{opp},E_n^{\Delta}) \cong HH_\idot(A,E_n) \cong
A^{\otimes n}.
$$
Since the category $\Lambda_{[n]}$ has an initial object $[n] \in
\Lambda_{[n]}$, we have $k = i_{n!k}$, so that the second multiple
$H_\idot(\Lambda_{[n]},k)$ is just $k$ in degree $0$. 
\endproof

The essential point of Proposition~\ref{main} is the following: the
cyclic object $A_\#$ associated to an algebra $A$ inconveniently
contains two things at the same time -- the cyclic structure, which
seems to be essential to the problem, and the bar resolution, which
is needed only to compute the Hochschild homology
$HH_\idot(A)$. Replacing $A_\#$ with the cyclic complex
$L^\hdot\tr_\#A_\# \in \D(\Lambda,k)$ disentagles these two. 

\medskip

We note that while one still has to prove that this does not change
the final answer, the construction itself looks pretty
straightforward -- if one wants to remove the non-essential bar
resolution from the definition of the cyclic homology,
Definition~\ref{cycl.cat} seems to be the obvious thing to
try. However, it was actually arrived at by a sort of a reverse
engeneering process. To finish the section, perhaps it would be
useful to show the reader the first stage of this process.

\medskip

Assume given an associative algebra $A$, and fix a projective
resolution $P_\idot$ of the diagonal $A$-module $A$. Then
$HH_\idot(A,M)$ can be computed by the complex
$$
\tr(P_\idot) = P_\idot \otimes_{A^{opp} \otimes A} A.
$$
How can one see the cyclic homology in terms of this complex? Or
even simpler -- what is the first differential in the spectral
sequence \eqref{conn.sp}, the Connes' differential $B:HH_\idot(A)
\to HH_{\idot+1}(A)$?

\medskip

There is the following recepy which gives the answer. Let
$\tau:P_\idot \to A$ be the augmentation map. Consider the tensor
product $P_\idot \otimes_A P_\idot$. This is also a projective
resotuion of $A$, and we actually have {\em two} natural
quasiisomorphisms
$$
\tau_1,\tau_2:P_\idot \otimes_A P_\idot \to P_\idot,
$$
given by $\tau_1 = \tau \otimes \id$, $\tau_2 = \id \otimes
\tau$. These quasiisomorphisms are different. However, since both
are maps between projective resolutions of the same object, there
should be a chain homotopy between them. Fix such a homotopy
$\iota:P_\idot \otimes_A P_\idot \to P_{\idot+1}$.

Now we apply the trace functor $\tr$, and obtain two maps
$\tau_1,\tau_2:\tr(P_\idot \otimes P_\idot) \to \tr(P_\idot)$, and a
homotopy $\iota:\tr(P_\idot \otimes P_\idot) \to \tr(P_{\idot+1})$
between them.

However, by the trace property of $\tau$, we also have an
involution $\sigma:\tr(P_\idot \otimes_A P_\idot)$ which
interchanges the two multiples. This involution obviously also
interchages $\tau_1$ and $\tau_2$, but there is no reason why it
should fix the homotopy $\iota$ -- in fact, it sends $\iota$ to a
second homotopy $\iota':\tr(P_\idot \otimes_A P_\idot) \to
\tr(P_{\idot+1})$ between $\tau_1$ and $\tau_2$.

The difference $\iota' - \iota$ is then a well-defined map of
complexes
\begin{equation}\label{homo}
\iota'-\iota:\tr(P_\idot \otimes_A P_\idot) \to \tr(P_{\idot+1}).
\end{equation}
On the level of homology, both sides are $HH_\idot(A)$; the map
$\iota'-\iota$ then induces exactly the Connes' differential
$B:HH_\idot(A) \to HH_{\idot+1}(A)$.

\medskip

To justify this recepy, we use Proposition~\ref{main} and identify
$HC_\idot(A)$ with $HC_\idot(L^\hdot\tr_\#(A_\#))$ rather than
$HC_\idot(A_\#)$. Then $L^\hdot\tr_\#(A_\#)$ is an object in
$\D_{const}(\Lambda,k)$. Therefore, as noted in
Remark~\ref{const.rem}, the Connes' differential $B$ only depends on
the restriction of $L^\hdot\tr_\#(A_\#)$ to $\Lambda_{\leq 2}
\subset \Lambda$. In other words, we do not need to compute the full
$L^\hdot\tr_\#(A_\#)$ and to construct a full resolution
$P^\#_\idot$ of the cyclic $A$-bimodule $A_\#$; it suffices to
construct $P^i_\idot = P^\#_\idot([i])$ for $i=1,2$ (and then apply
the functor $\tr$).

With the choices made above, we set $P^1_\idot = P_\idot$, and we
let $P^2_\idot$ be the cone of the map
$$
\begin{CD}
P_\idot \boxtimes P_\idot @>{(\tau \boxtimes \id)\oplus(\id
  \boxtimes \tau)}>> (A \boxtimes P_\idot)\oplus(P_\idot \boxtimes
  A).
\end{CD}
$$
The involution $\sigma:[2] \to [2]$ acts on $P^2_\idot$ in the
obvious way. We also need to define the transition maps $\iota_f$
for the two injections $d,d':[1] \to [2]$ and the two surjections
$s,s':[2] \to [1]$. For $d_1$, the transition map $\iota_d:A
\boxtimes P_\idot \to P^2_\idot$ is the obvious embedding, and so is
the transition map $\iota_{d'}$. For the surjection $s$, we need a
map $\iota_s$ from the cone of the map
$$
\begin{CD}
P_\idot \otimes_A P_\idot @>{(\tau \otimes \id)\oplus(\id
  \otimes \tau)}>> P_\idot \oplus P_\idot
\end{CD}
$$
to $P_\idot$. On $P_\idot \oplus P_\idot$, the map $\iota_s$ is just
the difference map $a \oplus b \mapsto a - b$; on
$P_\idot\otimes_AP_\idot$, $\iota_s$ is our fixed homotopy
$\iota:P_\idot\otimes_AP_\idot \to P_{\idot+1}$. And similarly for
the other surjection $s'$.

We leave it to the reader to check that if one computes
$L^\hdot\tr_\#(A_\#)\mid_{\Lambda_{\leq 2}}$ using this resolution
$P^\#_\idot$, then one obtains exactly \eqref{homo} for the Connes'
differential $B$.

\section{Discussion}

One of the most unpleasant features of the construction presented in
Section~\ref{cat} is the strong assumptions we need to impose on the
tensor category $\C$. In fact, the category one would really like to
apply the construction to is the category $\End\B$ of endofunctors
-- whatever that means -- of the category $\B$ of coherent sheaves
on an algebraic variety $X$. But if $X$ is not affine, $\End\B$
certainly does not have enough projectives, so that
Example~\ref{cln.exa} does not apply, and it is unlikely that
$\End\B$ can be made homologically clean in the sense of
Definition~\ref{clean}. We note that Definition~\ref{clean} has been
arranged so as not impose anything more than strictly necessary for
the proofs; but in practice, we do not know any examples which are
not covered by Example~\ref{cln.exa}.

As for the category $\End\B$, there is an even bigger problem with
it: while there are ways to define endofunctors so that $\End\B$ is
an abelian category with a right-exact tensor product, it cannot be
equipped with a right-exact trace functor $\tr$. Indeed, it
immediately follows from Definition~\ref{cycl.cat} that the
Hochschild homology groups $HH_\idot(\C)$ of a tensor category $\C$
are trivial in negative homological degrees. If $\C = \End\B$, one
of course expects $HH_\idot(\C) = HH_\idot(X)$, the Hochschild
homology $HH_\idot(X)$ of the variety $X$, which by now is
well-understood (see e.g. \cite{w}). And if $X$ is not affine,
$HH_\idot(X)$ typically is non-trivial both in positive and in
negative degrees. If $X$ is smooth and proper, $HH_\idot(X)$ in fact
carries a non-degenerate pairing, so that it is just as non-trivial
in degrees $>0$ as in degrees $<0$. Thus the case of a non-affine
algebraic variety is far beyond the methods developed in this paper.

The real reason for these difficulties is that we are dealing with
abelian categories, while the theory emphatically wants to live in
the triangulated world; as we explained in Example~\ref{const.exa},
even our main topic, cyclic bimodules, are best understood as
objects of a triangulated category $\D\Lambda(\C)$.  Unfortunately,
we cannot develop the theory from scratch in the triangulated
context, since we do not have a strong and natural enough notion of
an enhanced triangulated category (and working with the usual
triangulated categories is out of the question because, for
instance, the category of triangulated functors between triangulated
categories is usually not a triangulated category itself). A
well-developed theory would probably require a certain compromise
between the abelian and the triangulated approach. We will return to
it elsewhere.

Another thing which is very conspicously not done in the present
paper is the combination of Section~\ref{cat} and
Section~\ref{getz}. Indeed, in Section~\ref{getz}, we are dealing
with cyclic homology in the straightforward naive way of
Section~\ref{naive}, and while we define the cyclic object
$\wh{A_\#}$ associated to a square-zero extension $\wt{A}$, we make
no attempt to find an appropriate category $\wh{\Sec(A\bimod_\#)}$
where it should live. This is essentially the reason why we cannot
go further than square-zero extensions. At present, sadly, we do not
really understand this hypothetical category
$\wh{\Sec(A\bimod_\#)}$.

One suspects that treating this properly would require studying
deformations in a much more general context -- instead of
considering square-zero extensions of an algebra, we should look at
the deformations of the abelian category of its modules, or at the
deformations of the tensor category of its bimodules. This brings us
to another topic completely untouched in the paper: the Hochschild
cohomology $HH^\hdot(A)$.

Merely {\em defining} Hochschild cohomology for an arbitrary tensor
category $\C$ is in fact much simpler than the definition of
$HH_\idot(\C)$, and one does not need a trace functor for this -- we
just set $HH^\hdot(\C)=\Ext^\hdot(\I,\I)$, where $\I \in \C$ is the
unit object. However, it is well understood by now that just as
Hochschild homology always comes equipped with the Connes'
differential, the spectral sequence \eqref{conn.sp}, and the whole
cyclic homology package, Hochschild cohomology should be considered
not as an algebra but as the so-called {\em Gerstenhaber} algebra;
in fact, the pair $HH_\idot(-),HH^\hdot(-)$ should form a version of
``non-commutative calculus'', as proposed for instance in
\cite{TT}. Deformations of the tensor category $\C$ should be
controlled by $HH^\hdot(\C)$, and the behaviour of $HH_\idot(\C)$
and $HC_\idot(\C)$ under these deformations reflects various natural
actions of $HH^\hdot(-)$ on $HH_\idot(-)$.

We believe that a convenient development of the ``non-commutative
calculus'' for a tensor category $\C$ might be possible along the
same lines as our Section~\ref{cat}. Just as our category
$\D\Lambda(\C)$ is defined as the category of sections of the
cofibration $\C_\#/\Lambda$, whose definition imitates the usual
cyclic object $A_\#$, one can construct a cofibration $\C^\#/\Delta$
which imitates the standard cosimplicial object computing
$HH^\hdot(A)$ -- for any $[n] \in \Delta$, $\C^\#([n])$ is the
category of polylinear right-exact functors from $\C^{n-1}$ to $\C$,
and the transition functors between various $\C^\#([n])$ are induced
by the tensor product on $\C$. Then one can define a triangulated
category $\D\Delta(\C)$, the subcategory in $\D(\Sec(\C^\#))$ of
complexes with cocartesian homology; the higher structures on
$HH^\hdot(\C)$ should be encoded in the structure of the category
$\D\Delta(\C)$, and relations between $HH_\idot(\C)$ and
$HH^\hdot(\C)$ should be reflected in a relation between
$\D\Lambda(\C)$ and $\D\Delta(\C)$. We will proceed in this
direction elsewhere. At present, the best we can do is to make the
following hopeful observation:
\begin{itemize}
\item the category $\Sec_{cart}(\C^\#)$ is naturally a {\em braided}
  tensor category over $k$.
\end{itemize}
The reason for this is very simple: if one writes out explicitly the
definition of $\Sec_{cart}(\C^{\#})$ along the lines of
Lemma~\ref{cycl.str}, one finds out that it coincides on the nose
with the Drinfeld double of the tensor category $\C$.

\bigskip

\noindent
{\sc Steklov Math Institute\\
Moscow, USSR}

\bigskip

\noindent
{\em E-mail address\/}: {\tt kaledin@mccme.ru}

\end{document}